\documentclass[twocolumn]{article}
\usepackage{amssymb}
\newcommand{\vc}[1]{\nca{#1}}

\newcommand{\ica}[1]{\mbox{\it #1}}
\newcommand{\rca}[1]{\mbox{#1}}
\newcommand{\nca}[1]{\mbox{\bf #1}}

\newenvironment{pagi}[1]{\begin{minipage}[t]{#1em}}{\end{minipage}}

\newtheorem{prb}{Problem}[section]

\newcommand{\N}{\mathbb{N}}

\newcommand{\R}{\mathbb{R}}
\newcommand{\Z}{\mathbb{Z}}

\newcommand{\var}{\ica{Var}\,}

\title{
A straightforward local-search optimization algorithm on the symmetric group
}
\author{Guillermo Morales-Luna \\
Computer Science, CINVESTAV-IPN \\
Mexico City \\
{\tt gmorales@cs.cinvestav.mx}}
\date{\today}
\begin{document}
\maketitle

\begin{abstract}
Given a real objective function defined over the symmetric group, a direct local-search algorithm is proposed, and its complexity is estimated. In particular for an $n$-dimensional unit vector we are interested in the permutation isometry that acts on this vector by mapping it into a cone of a given angle.
\end{abstract}

\section{Introduction}

In several areas of Game Theory and Oil Engineering there appears an optimization problem defined over the group of permutations on an index set. The declination curves~\cite{LH03} on wells, for instance, are important in order to develop  exploitation strategies, among them Steam Assisted Gravitational Drainage (SAGD) processes~\cite{QG02}. The well known theories developed by Von Neumann on Game Theory pose also optimization problems whose domain is the symmetric group. A rather complete survey of this kind of optimization procedures in Social Sciences is given in~\cite{MM06}. 

A set of independent measurements of a certain parameter is realized as a vector in $\R^n$ and when it is normalized (in the statistical sense: it is translated by its average and elongated by the reciprocal of its variance) an unit vector with zero mean and unit variance is obtained. Then the problem consists in finding a permutation of its entries making a given angle with the original normalized vector.

Here we sketch a local-search procedure, with cubic complexity with respect to the dimension of the original vector.

\section{Uniform distribution}
Let $\vc{e}_i=\left(\delta_{ij}\right)_{j=1}^n$ be the $i$-th vector in the canonical basis of $\R^n$, where $\delta_{ij}$ is {\em Kroenecker's delta}. 
For any non-zero vector $\vc{x}=(x_1,\ldots,x_n)=\sum_{j=1}^n x_j\,\vc{e}_j$, its {\em mean} is $E(\vc{x}) = \frac{1}{n}\sum_{j=1}^n x_j$, its {\em variance} is $\var(\vc{x})= \frac{1}{n}\sum_{j=1}^n \left(x_j - E(\vc{x})\right)^2$ and its {\em normalization} is 
\begin{equation}
\vc{z}_{\bf x} = \frac{1}{\sqrt{\var(\vc{x})}}(\vc{x}-E(\vc{x})\,1^n) \label{eq.12}
\end{equation}
where $1^n$ is the vector in $\R^n$ with constant components 1. Clearly, $E(\vc{z}_{\bf x}) = 0$, $\var(\vc{z}_{\bf x}) = 1$ and $\vc{z}_{\bf x}$ is an unit vector, with respect to the Euclidean norm. In particular, for the {\em identity} vector
\begin{equation}
\vc{x} = (1,2,\ldots,n) = \sum_{j=1}^n j\,\vc{e}_j. \label{eq.01}
\end{equation}
its mean and variance are
\begin{eqnarray}
E(\vc{x}) &=& \frac{1}{n}\sum_{j=1}^n j = \frac{n+1}{2} \label{eq.02} \\
\var(\vc{x})&=& \frac{1}{n}\sum_{j=1}^n \left(j - E(\vc{x})\right)^2 = \frac{(n-1)(n+1)}{12}, \label{eq.021}
\end{eqnarray}
and its normalization is $\vc{z}=\sum_{j=1}^n z_j\vc{e}_j$, where, for each $j\leq n$:
\begin{equation}
z_j = \frac{1}{\sqrt{n}}\cdot \frac{j-E(\vc{x})}{\sqrt{\var(\vc{x})}} = \sqrt{\frac{3}{(n-1)n(n+1)}} ( 2 j - (n+1) ).\label{eq.03}
\end{equation}
Let
\begin{equation}
\vc{u} = \sqrt{\frac{(n-1)n(n+1)}{3}} \vc{z} = \sum_{j=1}^n \left( 2 j - (n+1) \right)\vc{e}_j. \label{eq.04}
\end{equation}
Then 
\begin{equation}
\|\vc{u}\|_2 = \sqrt{\frac{(n-1)n(n+1)}{3}}. \label{eq.05}
\end{equation}

\section{Optimization problem}
Let $S_n$ be the permutation group of the set of indexes $[\![1,n]\!]$. Then $\rca{card}\,S_n=n!$. 
A {\em permutation matrix} $P\in \R^{n\times n}$ has the form $P=\left(\delta_{i\pi(j)}\right)_{i,j\leq n}$ for some permutation $\pi\in S_n$. In this case, let us write $P=P_{\pi}$. This matrix determines a linear isometry $P_{\pi}:\R^n\to\R^n$: for each $\vc{y}\in\R^n$, $P_{\pi}(\vc{y})$ coincides with $\vc{y}$ up to the permutation $\pi$, hence $\|P_{\pi}(\vc{y})\|_2 = \|\vc{y}\|_2$. 

For any non-zero vector $\vc{x}=\sum_{j=1}^n x_j\,\vc{e}_j\in\R^n$, by Schwartz's inequality, the inner product between the normalized vector and its corresponding image satisfies
$\left|\langle \vc{z}_{\bf x}, P_{\pi}(\vc{z}_{\bf x})\right| \leq \|\vc{z}_{\bf x}\|_2\,\|P_{\pi}(\vc{z}_{\bf x})\|_2 = 1.$
Thus, for all $\pi\in S_n$,
\begin{equation}
-1\leq \langle \vc{z}_{\bf x}, P_{\pi}(\vc{z}_{\bf x}) \rangle \leq 1 \label{eq.07}
\end{equation}
Let us consider the map
\begin{equation}
\Phi_{\bf x}: S_n\to [-1,1]\ ,\ \pi \mapsto \langle \vc{z}_{\bf x}, P_{\pi}(\vc{z}_{\bf x}) \rangle \label{eq.08}
\end{equation}
For the {\em identity} permutation, $\pi_1=(1\ 2\ \cdots\ n)$ we have $\Phi_{\bf x}(\pi_1)=\langle \vc{z}_{\bf x}, \vc{z}_{\bf x} \rangle=1$ and for the permutation $\pi_0\in S_n$ such that $P_{\pi_0}(\vc{z}_{\bf x})=-\vc{z}_{\bf x}$ we have $\Phi_{\bf x}(\pi_0)=-1$.

In particular, if $\vc{x}$ is the identity vector, $\pi_0$ above is the {\em reverse} permutation $\pi_{n!}=(n\ n-1\ \cdots\ 1)$, namely $\forall j$, $\left( 2 \pi_{n!}(j) - (n+1) \right) = -\left( 2 j - (n+1) \right)$.
Besides, from eq's~(\ref{eq.03}) and~(\ref{eq.04}), 
\begin{eqnarray}
\langle \vc{z}, P_{\pi}(\vc{z}) \rangle &=& \frac{3}{(n-1)n(n+1)}\cdot\nonumber \\ 
&& \sum_{j=1}^n \left( 2 j - (n+1) \right) \left( 2 \pi(j) - (n+1) \right) \nonumber \\ &=& \frac{3}{(n-1)n(n+1)} \langle \vc{u}, P_{\pi}(\vc{u}) \label{eq.06}
\end{eqnarray}
Let us write $\nu_n=\frac{(n-1)n(n+1)}{3}$. 
From eq.~(\ref{eq.06}) it can be checked that the image of $\Phi$ is included in the discrete set
\begin{equation}
R_n =\left\{\left. \frac{k}{\nu_n}\right| k\in\Z\ ,\ -\nu_n \leq k \leq \nu_n\right\} \subset [-1,1] \label{eq.09}
\end{equation}
Namely, $\Phi(S_n)\subseteq R_n$, although the inclusion may be proper, for some values of $n$. For instance, by eq.~(\ref{eq.06}), $\forall k\in\Z$,
\begin{eqnarray}
&&\frac{k}{\nu_n}\in \Phi(S_n)  \Longleftrightarrow\ \nonumber \\ 
&&\exists\pi\in S_n: k = \sum_{j=1}^n \left( 2 j - (n+1) \right) \left( 2 \pi(j) - (n+1) \right) \label{eq.10}
\end{eqnarray}
Whenever $n$ is an even integer, each summand $\left( 2 j - (n+1) \right) \left( 2 \pi(j) - (n+1) \right)$ is odd (as the product of two odd integers), thus their addition will be even. Hence, for $k$ odd, we have $\frac{k}{\nu_n}\not\in \Phi(S_n)$.

Besides, we can see that the image $\Phi(S_n)$ is symmetric with respect to 0, namely, $\forall \pi\in S_n$
\begin{eqnarray}
k=\sum_{j=1}^n \left( 2 j - (n+1) \right) \left( 2 \pi(j) - (n+1) \right)&\Longrightarrow& \nonumber \\ 
-k=\sum_{j=1}^n \left( 2 j - (n+1) \right) \left( 2 \pi_{n!}\circ\pi(j) - (n+1) \right).&& \label{eq.11}
\end{eqnarray}

Let us consider the following
\begin{prb}\label{pb.01} Given a non-zero $\vc{x}\in\R^n$ and $r\in[-1,1]$ find $$\Pi(\vc{x},r)=\rca{\rm Arg Min}_{\pi\in S_n} |r-\langle \vc{z}_{\bf x}, P_{\pi}(\vc{z}_{\bf x},) \rangle |.$$
\end{prb}
From remark~(\ref{eq.11}), it is enough to assume $r\in[0,1]$ when $\vc{x}$ is the identity vector.

Problem~\ref{pb.01} can be solved by a direct {\em local-search} algorithm. Let us consider the elements of $S_n$ as {\em nodes} and let us consider {\em edges} of two kinds:
\begin{eqnarray*}
(\pi,\rho)\in E_T &\Longleftrightarrow& \exists i,j\in\{1,\ldots,n\}:\ \rho =(i,j)\pi \\
(\pi,\rho)\in E_C &\Longleftrightarrow& \exists i\in\{1,\ldots,n\}:\ \rho =(i,i')\pi \ \nonumber \\ 
&&\mbox{with }i'=i+1\mbox{ if }i<n\ \&\ n'=1
\end{eqnarray*}
Let ${\bf S}_{nT}=(S_n,E_T)$ and ${\bf S}_{nC}=(S_n,E_C)$. The graph ${\bf S}_{nT}$ contains as edges all pairs of permutations that differ by an arbitrary transposition, and the graph ${\bf S}_{nC}$ contains as edges all pairs of permutations that differ by a transposition of two consecutive indexes. ${\bf S}_{nT}$ is a regular graph of $n!$ nodes of degree ${n \choose 2} =\frac{1}{2}n(n-1)$. If we consider any edge as an edge of unit length, then the {\em diameter} of ${\bf S}_{nT}$ is $d_{nT}=n-1=O(n)$. Instead, ${\bf S}_{nC}$ is a regular graph of $n!$ nodes of degree $n$ and its diameter is 
$$d_{nC}=\frac{1}{4}(n-(n\,\mbox{mod}\,2))(n+(n\,\mbox{mod}\,2))=O(n^2).$$ 

\section{Solution procedure}

Let $F:S_n\to\R^+$ be a positive-real valued map defined over the symmetric group. 
Let ${\bf S}_{n}=(S_n,E)$ be any of ${\bf S}_{nT}$ or ${\bf S}_{nC}$. A local-search algorithm to minimize $F$ is sketched at table~\ref{tb.01}

\begin{table}
\begin{center}
 \fbox{\begin{minipage}{24em}
 \noindent{\bf Input.} The dimension $n\in\N$ and the objective map $F:S_n\to\R^+$.

\noindent{\bf Output.} The permutation $\pi\in S_n$ such that $\pi=\rca{\rm Arg Min}_{\rho\in S_n} F(\rho)$.

\begin{enumerate}
\item  Initialize $\pi_c:=\pi_1$ (the identity permutation); $v_c:=F(\pi_c)$;
\item  repeat
\begin{enumerate}
\item  $\pi_d:= \pi_c$; $v_d:=v_c$;
\item  let $\pi_c$ be the $E$-neighbor of $\pi_d$ minimizing the objective map $F$ in the $E$-neighborhood of $\pi_d$;
\item  $v_c:=F(\pi_c)$;
\item  if $v_c\geq v_d$ then $\pi_c:= \pi_d$;
\end{enumerate}
until $\pi_c==\pi_d$;
\item  output $\pi:=\pi_c$.
\end{enumerate}
\end{minipage}}
\end{center}
 \caption{\label{tb.01}Local-search algorithm on the symmetric group.}
\end{table}

In order to solve problem~\ref{pb.01}, for any given $\vc{x}\in\R^n-\{\vc{0}\}$, $r\in[0,1]$, the objective map $F_{{\bf x},r}:\pi\mapsto |r- \Phi_{\bf x}(\pi)|$ is considered, where $\Phi_{\bf x}$ is given by eq.~(\ref{eq.08}) which in turn involves the normalization $\vc{z}_{\bf x}$ of vector $\vc{x}$, given by eq.~(\ref{eq.03}) for the identity vector.

If ${\bf S}_{n}={\bf S}_{nT}$, then the graph is regular of degree $n(n-1)/2$ and diameter $n-1$. Thus step 2.(b) at pseudocode in table~\ref{tb.01} entails $n(n-1)/2$ evaluations of the objective map, and whole cycle 2. is repeated at most $n-1$ times. Thus, in this case the worst case entails $n(n-1)^2/2$ evaluations of the objective map.

Instead, if ${\bf S}_{n}={\bf S}_{nC}$, then the graph is regular of degree $n$ and diameter $d_{nC}$. Thus step 2.(b) at pseudocode in table~\ref{tb.01} entails $n$ evaluations of the objective map, and whole cycle 2. is repeated at most $d_{nC}$ times. Thus, in this case the worst case entails $nd_{nC}$ evaluations of the objective map. Since
$$\frac{n(n-1)^2/2}{n(n-(n\,\mbox{mod}\,2))(n+(n\,\mbox{mod}\,2))/4} = 2\frac{(n-1)^2}{(n-(n\,\mbox{mod}\,2))(n+(n\,\mbox{mod}\,2))}\mathop{\longrightarrow}_{n\to+\infty} 2,$$
we see that using ${\bf S}_{nT}$ is twice as expensive than using ${\bf S}_{nC}$, and this last entails $O(n^3)$ evaluations of the objective map.

\end{document}